\documentclass[a4paper,reqno]{amsart}
\usepackage[all]{xy} %For commutative diagrams
\usepackage{amssymb} %For double arrows
\usepackage{hyperref}
\usepackage{eucal}
\usepackage{pb-diagram}
\usepackage{graphicx}
\usepackage{epsfig}
\numberwithin{equation}{section}

\newtheorem{definition}{Definition}[section]

\newtheorem{theorem}[definition]{Theorem}

\newtheorem{remarkth}[definition]{Remark}
\newtheorem{example}[definition]{Example}

\renewcommand{\emph}[1]{{\bfseries\itshape{#1}}}

\newcommand{\R}{\mathbb{R}} %Numeros reales
\makeatletter
\newcommand\prol{\@ifstar{\@proldf}{\@prolpf}} %% if * dual else primal
\def\@prolpf{\@ifnextchar[{\@prolpf@wrt}{\@prolpf@}}
\def\@prolpf@wrt[#1]#2{\@ifnextchar[{\@prolpf@wrt@at{#1}{#2}}{\@prolpf@wrt@{#1}{#2}}}

\def\@prolpf@wrt@at#1#2[#3]{\prolsymbol^{#1}_{#3}#2}
\def\@prolpf@wrt@#1#2{\prolsymbol^{#1}#2}
\def\@prolpf@#1{\@ifnextchar[{\@prolpf@at{#1}}{\@prolpf@@{#1}}}
\def\@prolpf@at#1[#2]{\prolsymbol_{#2}#1}
\def\@prolpf@@#1{\prolsymbol#1}
\def\@proldf{\@ifnextchar[{\@proldf@wrt}{\@proldf@}}
\def\@proldf@wrt[#1]#2{\@ifnextchar[{\@proldf@wrt@at{#1}{#2}}{\@proldf@wrt@{#1}{#2}}}

\def\@proldf@wrt@at#1#2[#3]{\prolsymbol^{*#1}_{#3}#2}
\def\@proldf@wrt@#1#2{\prolsymbol^{*#1}#2}
\def\@proldf@#1{\@ifnextchar[{\@proldf@at{#1}}{\@proldf@@{#1}}}
\def\@proldf@at#1[#2]{\prolsymbol^*_{#2}#1}
\def\@proldf@@#1{\prolsymbol^*#1}
\def\prolsymbol{\mathcal{L}}
\makeatother

\newcommand{\Sec}[1]{\operatorname{Sec}(#1)}

\newcommand{\id}{{\operatorname{id}}}
\font\blackital=cmmib10  \skewchar\blackital='177
\font\sblackital=cmmib7  \skewchar\sblackital='177
\font\ssblackital=cmmib5  \skewchar\ssblackital='177
\font\sanss=cmss10 \font\ssanss=cmss8 scaled 900
\font\sssanss=cmss8 scaled 600 
 \font\teneurm=eurm10 
\font\blackboard=msbm10 \font\sblackboard=msbm7
\font\ssblackboard=msbm5 \font\caligr=eusm10 \font\scaligr=eusm7
\font\sscaligr=eusm5  \font\fraktur=eufm10
\font\sfraktur=eufm7 \font\ssfraktur=eufm5 
\newfam\bbfam
\textfont\bbfam=\blackboard \scriptfont\bbfam=\sblackboard
\scriptscriptfont\bbfam=\ssblackboard

\newfam\ssfam
\textfont\ssfam=\sanss \scriptfont\ssfam=\ssanss
\scriptscriptfont\ssfam=\sssanss
\def\ss#1{{\fam\ssfam\relax#1}}

\newfam\clfam
\textfont\clfam=\caligr \scriptfont\clfam=\scaligr
\scriptscriptfont\clfam=\sscaligr
\def\cl#1{{\fam\clfam\relax#1}}

\newfam\frfam
\textfont\frfam=\fraktur \scriptfont\frfam=\sfraktur
\scriptscriptfont\frfam=\ssfraktur

\newfam\gpfam
\def\eu#1{\hbox{$\fam\gpfam\relax#1\textfont\gpfam=\teneurm$}}

\font\bsymb=cmsy10 scaled\magstep2
\def\all#1{\setbox0=\hbox{\lower1.5pt\hbox{\bsymb \char"38}}\setbox1=\hbox{$_{#1}$} \box0\lower2pt\box1}
\def\exi#1{\setbox0=\hbox{\lower1.5pt\hbox{\bsymb \char"39}}\setbox1=\hbox{$_{#1}$} \box0\lower2pt\box1}

\catcode`\"=12
\mathchardef\za="710B  %\alpha
\mathchardef\zb="710C  %\beta
\mathchardef\zg="710D  %\gamma
\mathchardef\zd="710E  %\delta
\mathchardef\zve="710F %\epsilon
\mathchardef\zz="7110  %\zeta
\mathchardef\zh="7111  %\eta
\mathchardef\zvy="7112 %\theta
\mathchardef\zi="7113  %\iota
\mathchardef\zk="7114  %\kappa
\mathchardef\zl="7115  %\lambda
\mathchardef\zm="7116  %\mu
\mathchardef\zn="7117  %\nu
\mathchardef\zx="7118  %\xi
\mathchardef\zp="7119  %\pi
\mathchardef\zr="711A  %\rho
\mathchardef\zs="711B  %\sigma
\mathchardef\zt="711C  %\tau
\mathchardef\zu="711D  %\upsilon
\mathchardef\zvf="711E %\phi
\mathchardef\zq="711F  %\chi
\mathchardef\zc="7120  %\psi
\mathchardef\zw="7121  %\omega
\mathchardef\ze="7122  %\varepsilon
\mathchardef\zy="7123  %\vartheta
\mathchardef\zvp="7124  %\varpi
\mathchardef\zvr="7125 %\varrho
\mathchardef\zvs="7126 %\varsigma
\mathchardef\zf="7127  %\varphi
\mathchardef\zG="7000  %\Gamma
\mathchardef\zD="7001  %\Delta
\mathchardef\zY="7002  %\Theta
\mathchardef\zL="7003  %\Lambda
\mathchardef\zX="7004  %\Xi
\mathchardef\zP="7005  %\Pi
\mathchardef\zS="7006  %\Sigma
\mathchardef\zU="7007  %\Upsilon
\mathchardef\zF="7008  %\Phi
\mathchardef\zC="7009  %\Psi
\mathchardef\zW="700A  %\Omega
\catcode`\"=\active

\def\ezt{\eu\zt}

\def\sT{{\ss T}}
\def\cC{{\cl C}}

\def\cR{{\cl R}}

\def\sss#1{\hbox{\ssanss #1}}

\newcommand{\be}{\begin{equation}}
\newcommand{\ee}{\end{equation}}
\newcommand{\ra}{\rightarrow}

\newcommand{\bea}{\begin{eqnarray}}
\newcommand{\eea}{\end{eqnarray}}
\newcommand{\beas}{\begin{eqnarray*}}
\newcommand{\eeas}{\end{eqnarray*}}
\def\*{{\textstyle *}}
\def\xd{\operatorname{d}\!}

\def\bt{\mathbf{t}}
\def\dt{\xd_{\sss T}}
\def\vt{\textsf{v}_{\sss T}}
\def\rel{{-\!\!\!-\!\!\rhd}}

\begin{document}

\title[Nonholonomic Constraints: a New Viewpoint]{Nonholonomic
Constraints: a New Viewpoint}

\author[J. Grabowski]{J. Grabowski}
\address{Janusz Grabowski: Polish Academy of Sciences, Institute of Mathematics,
 \'Sniadeckich 8, P.O. Box 21, 00-956 Warszawa, Poland
} \email{jagrab@impan.gov.pl}

\author[M. de Le\'on]{M. de Le\'on}
\address{Manuel de Le\'on:
Instituto de Ciencias Matem\'aticas (CSIC-UAM-UC3M-UCM),  Serrano 123, 28006
Madrid, Spain} \email{mdeleon@imaff.cfmac.csic.es}

\author[J.\ C.\ Marrero]{J.\ C.\ Marrero}
\address{Juan Carlos Marrero:
Departamento de Matem\'atica Fundamental y Unidad Asociada
ULL-CSIC Geo\-metr\'{\i}a Diferencial y Mec\'anica Geom\'etrica,
Facultad de Ma\-te\-m\'a\-ti\-cas, Universidad de la Laguna, La
Laguna, Tenerife, Canary Islands, Spain} \email{jcmarrer@ull.es}

\author[D.\ Mart\'{\i}n de Diego]{D. Mart\'{\i}n de Diego}
\address{David Mart\'{\i}n de Diego:
Instituto de Ciencias Matem\'aticas (CSIC-UAM-UC3M-UCM),  Serrano 123, 28006
Madrid, Spain} \email{d.martin@imaff.cfmac.csic.es}

\keywords{General algebroids, Lie algebroids, double vector
bundle, Nonholonomic mechanics, Lagrange-d'Alembert's equations,
nonholonomic bracket.}

\subjclass[2000]{70F25, 53D17, 70G45, 17B66, 70H03, 70H45}

\date{}

\thanks{This work has been partially supported by the Polish
Ministry of Science and Higher Education
under the grant No. N201 005 31/0115, MEC (Spain)
Grants MTM 2006-03322, MTM 2007-62478, project Ingenio
Mathematica (i-MATH) No. CSD 2006-00032 (Consolider-Ingenio 2010)
and S-0505/ESP/0158 of the CAM}

\begin{abstract}
The purpose of this paper is to show that, at least for
Lagrangians of mechanical type, nonholonomic Euler-Lagrange
equations for a nonholonomic linear constraint $D$ may be viewed
as non-constrained Euler-Lagrange equations but on a new
(generally not Lie) algebroid structure on $D$. The proposed novel
formalism allows us to treat in a unified way a variety of
situations in nonholonomic mechanics and gives rise to a version
of Neoether Theorem producing actual first integrals in case of
symmetries.
\end{abstract}

\maketitle

\tableofcontents
\section{Introduction}
There are many approaches to geometric mechanics in the
literature. We will work with a natural generalization of the
framework for studying mechanical systems proposed by
W.~M.~Tulczyjew \cite{Tul1,Tul2} (see also \cite{TU} and
references therein). In the simplest form, the phase dynamics of
the system is understand as the lagrangian submanifold $\zG$ of
the symplectic manifold $(\sT\sT^*M,\dt\omega_M)$ equipped with
tangent lift $\dt\omega_M$ of the canonical symplectic form
$\omega_M$ of $\sT^*M$. Here, $M$ represents the configuration
manifold of the system and $\zG$ is obtained from the lagrangian
submanifold $\xd L(M)\subset \sT^*\sT M$ induced by the Lagrangian
$L:\sT M\ra\R$ via the canonical isomorphism $\ze_M:\sT^*\sT M\ra
\sT\sT^*M$. In other words, the phase dynamics, as well as the
Euler-Lagrange equations, are obtained in a simple way by means of
the {\it Tulczyjew differential} $\zL_L=\ze_M\circ\xd
L:M\ra\sT\sT^* M$. It is important to observe that both
$\sT\sT^*M$ and $\sT^*\sT M$ are double vector bundles over
$\sT^*M$ and $\sT M$ (see \cite{GU2} and references therein). The
resulting submanifold $\zG$ of $\sT\sT^*M$ is a particular case of
modelling dynamical systems as implicit differential equations
defined by differential inclusions (see \cite{MMT1,MMT2}).

This framework admits an immediate generalization for more general
morphisms $\ze:\sT^* E\ra\sT E^*$ of canonical double vector
bundles associated with a vector bundle $E$ (see \cite{GGU,GG}),
inducing the Tulczyjew differentials $\zL_L=\ze\circ\xd L:M\ra\sT
E^*$. This generalization includes as a particular case a theory
of mechanical systems based on Lie algebroids, as proposed by
A.~Weinstein \cite{We} and developed by many authors, however in a
different geometrical setting (see, for instance,
\cite{CoLeMaMaMa,CoLeMaMa,LMM,Li,mart,Medina}). The motivation for
study systems on Lie algebroids is that they often appear
naturally as results of some reduction procedures. This is is a
situation similar to the one known in the theory of Hamiltonian
systems: reductions may lead from a symplectic to a Poisson
structure.

An additional challenge and one of the most fascinating topics in
geometric mechanics is the study of constraints in this context.
Of course, a general problem of putting constraints for the system
in a variational setting involves constraints for velocities as
well as constraints for virtual displacements, as was noticed
already in \cite{Tu2}. In some cases, however, one assumes that
the constraints can be determined from a constraint subset $D$ of
$\sT M$ (or, of $E$ in the algebroid case) by certain
well-described procedures. The best known approaches of this type
refer to the so called {\it vakonomic} and {\it nonholonomic}
constraints. In the simplest situation, for $D$ being a linear
nonholonomic constraint, i.e. just a vector subbundle of $\sT M$
(or, of $E$ in the algebroid context), this procedure describes
the nonholonomic Euler-Lagrange equations by means of the {\it
d'Alembert principle}, having analogs also in the algebroid case
\cite{CoLeMaMa,LMMdD,M4,GG}. We should stress that our
nonholonomic constraints are linear in the broader sense, i.e.
they are subbundles over submanifolds of the original base
manifold. The nonholonomic Euler-Lagrange equations are commonly
viewed as being not variational equations. In \cite{GG} it has
been pointed out that it is not exactly the case, if we extend
slightly our understanding of Variational Calculus.

In this paper we continue studying nonholonomic constraints on
algebroids and showing that, at least for Lagrangians of
mechanical type, the nonholonomic Euler-Lagrange equations are
just non-constrained Euler-Lagrange equations but for a special
algebroid structure on the constraint subbundle (see also
\cite{LMMdD}). This shows that mechanical systems based on general
(not necessary Lie) algebroids appear naturally in the presence of
nonholonomic constraints and gives a powerful geometrical tool
when dealing with constrained systems. In particular, we get a
version of Noether Theorem with true first integrals for
nonholonomic systems. We do not get all possible algebroids on $D$
applying our procedure. In particular, if the original structure
was a Lie algebroid, then the new algebroid bracket is
automatically skew- symmetric, so we deal with a {\it quasi-Lie
algebroid}. One can the associate with the sequence of procedures,
like reduction by symmetries and passing to a nonholonomic
constraint, the sequence of the corresponding novel structures
serving as appropriate geometrical tools in describing the
systems:
$$ \sT M\quad \xrightarrow[]{\text{reduction by symmetries}}\quad \text{Lie algebroid} \quad
\xrightarrow[]{\text{nonholonomic constraint}}\quad\text{quasi-Lie
algebroid}\,.
$$
All this is of course closely related to the discovery of the role
of the nonholonomic quasi-Poisson brackets
\cite{Mr,SM,IbLeMaMa,CaLeMa}, this time not in the Hamilton but in
the Lagrange picture.

The paper is organized as follows. In the next section, we recall
after \cite{GGU} the basic ideas of developing mechanics on a
general algebroid $E$, in particular, the Euler-Lagrange
equations. Then, we construct in this setting an analog of the
Tulczyjew differential for linear nonholonomic constraints,
together with the corresponding nonholonomic Euler-Lagrange
equations. In Section 4 we study reductions of a general algebroid
to an algebroid on a nonholonomic constraint $D$ (satisfying a
natural admissibility condition) along a given projection. We
discuss also the problem, what algebroids can be obtained in this
way, if we start with Lie algebroids.

Section 5 is the most important part of our paper. For Lagrangian
functions $L$ of mechanical type (`kinetic energy - potential'),
we show that the non-constrained Euler-Lagrange equations on $D$,
derived for the reduced Lagrangian $l=L_{\mid D}$ and for the
reduced algebroid structure on $D$ along the orthogonal projection
associated with the kinetic energy, coincide with the nonholonomic
Euler-Lagrange equations. Moreover, for this nonholonomic case, we
can apply therefore the generalization of the Noether Theorem
proved in \cite{GGU} to obtain actual first integrals. Passing to
the nonholonomic constraints does not requires therefore any
change in our unified algebroid approach to mechanics:
nonholonomic Euler-Lagrange equations are included in our
framework.

We end up with two well-known examples of nonholonomic
constraints, the Chaplygin sleigh and the snakebord, to show how
simply the corresponding equations of motions can be derived by
means of our method.
\section{Geometric Mechanics on general algebroids}

Let $M$ be a smooth manifold and $(x^i)$, $i=1, \ldots, n$, be a
coordinate system in $M$. Denote by $\tau_M\colon \sT M\rightarrow
M$ the tangent vector bundle, with induced coordinates $(x^i,
\dot{x}^j)$, and by $\pi_M\colon \sT^\ast M \rightarrow M$ the
contangent bundle, with induced coordinates $(x^i, p_j)$.

Let $\tau: E\rightarrow M$ be a vector bundle and $\pi: E^\ast\rightarrow M$ its dual
bundle. Taking a local basis $\{e_a\}$ of sections of $E$, then we have the
corresponding local coordinates $(x^i, y^a)$ on $E$, where $y^a(e)$ is the
$a$th-coordinate of $e\in E$ in the given basis. We denote by $(x^i, \zx_a)$ the
corresponding coordinates of the dual bundle $\pi: E^\ast\rightarrow M$. One can also
say that $\zx_a$ is the fiber-wise linear local function $\zi(e_a)$ on $E^\ast$
corresponding to the local section $e_a$ of $E$. We have also adapted local
coordinates:
    \beas
     (x^i, y^a,{\dot x}^j, {\dot y}^b ) \quad \text{in} \ \sT E\,,& \qquad
    (x^i, \zx_a, {\dot x}^j, {\dot \zx}_b) \quad \text{in} \ \sT E^\* \,,\\
    (x^i, y^a, p_j, \zp_b)\quad \text{in}\ \sT^\*E\,,&\qquad
    (x^i, \zx_a, p_j, \zf^b)\quad \text{in}\ \sT^\* E^\* \,.
    \eeas

It is well known (cf. \cite{KU}) that the cotangent bundles $\sT^\*E$ and $\sT^\*E^\*$
are examples of double vector bundles:
$$\xymatrix{
\sT^\ast E^\ast\ar[rr]^{\sT^\ast\zp} \ar[d]_{\zp_{E^\ast}} && E\ar[d]^{\zt} \\
E^\ast\ar[rr]^{\zp} && M } \qquad {,}\qquad \xymatrix{
\sT^\ast E\ar[rr]^{\sT^\ast\zt} \ar[d]_{\zp_{E}} && E^\ast\ar[d]^{\zp} \\
E\ar[rr]^{\zt} && M }.
$$
The concept of a double vector bundle is due to J.~Pradines
\cite{Pr1,Pr2}, see also \cite{Ma,KU}. In particular, all arrows
correspond to vector bundle structures and all pairs of vertical
and horizontal arrows are vector bundle morphisms. The double
vector bundles have been recently characterized \cite{GR} in a
simple way as two vector bundle structures whose Euler vector
fields commute. The above double vector bundles are canonically
isomorphic with the isomorphism
    \be\label{iso}\cR_\zt \colon \sT^\*E \rightarrow \sT^\* E^\*
                    \ee
  being simultaneously an anti-symplectomorphism
  (cf. \cite{KU,GU2}). In adapted local coordinates, $\cR_\zt$ is given by
    $$\cR_\zt(x^i, y^a, p_j, \zp_b) = (x^i, \zp_a, -p_j,y^b).
                              $$
This means that we can identify $\zp_b$ and $\zx_b$, as well as $y^b$ with $\zf^b$, so
we will use local coordinates $(x,y,p,\zx)$ on $\sT^\ast E$ and local coordinates
$(x,\zx,p,y)$ on $\sT^\ast E^\ast$ in full agreement with the isomorphism (\ref{iso}).

For describing the dynamics of a Lagrangian or Hamiltonian systems
it is necessary to give an additional ingredient, typically, a Lie
algebroid structure or, more generally, the algebroid structure.
In several papers, many authors have studied Lie algebroid
structures as a unified geometric framework, general enough to
account for different mechanical systems (defined on tangent
bundles, on Lie algebras, quotients by Lie groups actions, etc.).

\medskip It is well known that Lie algebroid structures on a vector bundle $E$
correspond to linear Poisson tensors on $E^\*$. A 2-contravariant tensor $\zP$ on
$E^\*$ is called {\it linear} if the corresponding mapping $\widetilde{\zP} \colon
\sT^\* E^\* \rightarrow \sT E^\*$ induced by contraction is a morphism of double vector
bundles. This is the same as to say that the corresponding bracket of functions is
closed on (fiber-wise) linear functions. The commutative diagram
$$\xymatrix{
\sT^\ast E^\ast\ar[r]^{\widetilde\Pi}  & \sT E^\ast \\
\sT^\ast E\ar[u]_{\cR_\tau}\ar[ur]^{\ze} & },
$$
composed with (\ref{iso}), describes a one-to-one correspondence between linear
2-contravariant tensors $\zP$ on $E^\*$ and homomorphisms of double vector bundles (cf.
\cite{KU, GU2}) covering the identity on $E^\*$:

\be\xymatrix{
 & \sT^\ast E \ar[rrr]^{\varepsilon} \ar[dr]^{\pi_E}
 \ar[ddl]_{\sT^\ast\tau}
 & & & \sT E^\ast\ar[dr]^{\sT\pi}\ar[ddl]_/-20pt/{\tau_{E^\ast}}
 & \\
 & & E\ar[rrr]^/-20pt/{\zr}\ar[ddl]_/-20pt/{\tau}
 & & & \sT M \ar[ddl]_{\tau_M}\\
 E^\ast\ar[rrr]^/-20pt/{id}\ar[dr]^{\pi}
 & & & E^\ast\ar[dr]^{\pi} & &  \\
 & M\ar[rrr]^{id}& & & M &
}\label{F1.3} \ee In local coordinates, every  $\ze$ as above is
of the form
\be\label{F1.4} \ze(x^i,y^a,p_j,\zp_b) = (x^i, \zp_a,
\sum_c\zr^j_c(x)y^c, \sum_{a,c}{\mathcal C}^c_{ab}(x) y^a\zp_c +
\sum_i\zs^i_b(x) p_i) \ee
and it corresponds to the linear tensor
\be\label{tensor} \zP_\ze =\sum_{a,b,c}{\mathcal C}^c_{ab}(x)\zx_c
\partial _{\zx_a}\otimes
\partial _{\zx_b} + \sum_{i,b}\zr^i_b(x) \partial _{\zx_b} \otimes
\partial _{x^i} - \sum_{a,j}\zs^j_a(x)\partial _{x^j} \otimes
\partial _{\zx_a}\,.
\ee
In \cite{GU2} by {\it algebroids} we meant the morphisms
(\ref{F1.3}) of double vector bundles covering the identity on
$E^\*$, while {\it Lie algebroids} were those algebroids for which
the tensor $\zP_\ze$ is a Poisson tensor. We can consider the {\it
adjoint tensor} $\zP_\ze^+$, i.e. the 2-contravariant tensor
obtained from $\zP_\ze$ by transposition:
$$ \zP_\ze^{+} =\sum_{a,b,c}{\mathcal C}^c_{ba}(x)\zx_c \partial _{\zx_a}\otimes
\partial _{\zx_b} + \sum_{i,b}\zr^i_b(x) \partial _{x^i}\otimes\partial _{\zx_b}
 - \sum_{a,j}\zs^j_a(x) \partial _{\zx_a}\otimes\partial _{x^j}\,,
$$
and the corresponding {\it adjoint algebroid} structure $\ze^+$.
Algebroids $\ze$ corresponding to skew-symmetric $\zP_\ze$
(anti-symmetric brackets $[\cdot,\cdot]_\ze$), i.e. such that
$\ze^{+}=-\ze$, we will call {\it quasi-Lie algebroids}. The
relation to the canonical definition of Lie algebroid is given by
the following theorem (cf. { \cite{GU3, GU2}}).

\begin{theorem}
An algebroid structure $(E,\ze)$ can be equivalently defined as a bilinear bracket $[\cdot ,\cdot]_\ze $ on
sections of $\zt\colon E\rightarrow M$, together with vector bundle morphisms $\zr,\, \zs \colon E\rightarrow
\sT M$ (left and right anchors), such that
$$ [fX,gY]_\ze = f\zr(X)(g)Y -g\zs(Y)(f) X
+fg [X,Y]_\ze
$$
         for $f,g \in \cC^\infty (M)$, $X,Y\in\Sec{E}$.
The bracket and anchors are related to the 2-contravariant tensor $\zP_\ze$ by the
formulae
\beas
        \zi([X,Y]_\ze)&= \{\zi(X), \zi(Y)\}_{\zP_\ze},  \\
        \zp^\*(\zr(X)(f))       &= \{\zi(X), \zp^\*f\}_{\zP_\ze}, \\
        \zp^\*(\zs(X)(f))       &= \{\zp^\* f, \zi(X)\}_{\zP_\ze}.
                                                   \eeas
        The algebroid $(E,\ze)$ is a quasi-Lie algebroid if and only if the tensor $\zP_\ze$ is
        skew-symmetric and it is a Lie algebroid if and only if the tensor
$\zP_\ze$ is a Poisson tensor.
\end{theorem}

The canonical example of a mapping $\ze$ in the case of $E=\sT M$
is given by $\ze = \ze_M = \za^{-1}_M$ -- the inverse to the
Tulczyjew isomorphism $\za_M:\sT\sT^\*M\ra\sT^\*\sT M$. In
general, the algebroid structure map $\ze$ is not an isomorphism
and, consequently, its dual $\zk = \ze^{\ast}$ with respect to the
right projection is a relation and not a mapping. Some authors
(see \cite{LMM} for the case of Lie algebroids) have introduced
the concept of prolongation of a Lie algebroid to maintain some of
the original properties of the Tulzcyjew mapping  (isomorphism and
symplecticity), but, for the purposes of this paper, it is not
necessary to use this formalism.

The double vector bundle morphism (\ref{F1.3}) can be extended to the following
algebroid analogue of the so called Tulczyjew triple
\be\xymatrix@C-5pt{
 &\sT^\ast E^\ast \ar[rrr]^{\widetilde{\zP}}
\ar[ddl]_{\pi_{E^\ast}} \ar[dr]^{\sT^\ast\pi}
 &  &  & \sT E^\ast \ar[ddl]_/-25pt/{\zt_{E^\ast}} \ar[dr]^{\sT\pi}
 &  &  & \sT^\ast E \ar[ddl]_/-25pt/{\sT^\ast\zt} \ar[dr]^{\pi_E}
\ar[lll]_{\varepsilon}
 & \\
 & & E \ar[rrr]^/-20pt/{\zr}\ar[ddl]_/-20pt/{\zt}
 & & & \sT M\ar[ddl]_/-20pt/{\zt_M}
 & & & E\ar[lll]_/+20pt/{\zr}\ar[ddl]_{\zt}
 \\
 E^\ast\ar[rrr]^/-20pt/{id} \ar[dr]^{\pi}
 & & & E^\ast\ar[dr]^{\pi}
 & & & E^\ast\ar[dr]^{\pi}\ar[lll]_/-20pt/{id}
 & & \\
 & M\ar[rrr]^{id}
 & & & M & & & M\ar[lll]_{id} &
}\label{F1.3b}\ee The left-hand side is Hamiltonian, the right-hand side is Lagrangian,
and the dynamics lives in the middle.
%tu

We introduce now the dynamics through a  Lagrangian $L:E\ra{\R}$
which defines two smooth maps: the {\em Legendre mapping}: {
$\lambda_L:E\rightarrow E^\ast$,
$\lambda_L=\tau_{E^\ast}\circ\ze\circ\xd L$,} which is covered by
the {\em Tulczyjew differential} { $\Lambda_{L}: E\rightarrow \sT
E^\ast$, $\Lambda_{L}=\ze\circ\xd L$}:
$$\xymatrix{
\sT^\ast E\ar[rr]^{\ze}  && \sT E^\ast\ar[d]^{\zt_{E^\ast}} \\
E\ar[rr]^{\zl_L}\ar[u]^{\xd L}\ar[rru]^{\Lambda_{L}} && E^\ast }.
$$
The lagrangian function defines {\em the phase dynamics} {
${\zG_L}=\Lambda_{L}(E)\subset \sT E^\ast$} which can be
understood as an implicit differential equation on $E^\*$,
solutions of which are `phase trajectories' of the system
$\zb:\R\ra E^\ast$ and satisfy $\bt{\zb}(t)\in {\zG_L}$, where
$\bt{\zb}$ is the tangent prolongation of $\zb$,
$\bt\zb(t)=(\zb(t),\dot{\zb}(t))$. An analog of the {\it
Euler-Lagrange equations} for curves $\gamma: \R\rightarrow E$ is
in turn:
$$(E_L)\qquad \bt(\zl_L\circ\gamma)=\Lambda_L\circ\gamma.$$
In local coordinates, ${\zG_L}$ has the parametrization by
$(x^i,y^a)$ via $\zL_L$ in the form (cf. (\ref{F1.4}))
\be\zL_L(x^i,y^a)= (x^a,\frac{\partial L}{\partial
y^a}(x,y), \sum_c\zr^j_c(x)y^c, \sum_{a,c}{\mathcal C}^c_{ab}(x) y^a\frac{\partial L}{\partial y^c}(x,y) +\sum_i
\zs^i_b(x)\frac{\partial L}{\partial x^i}(x,y)) \label{F1.4a}\ee and the equation $(E_L)$, for
$\zg(t)=(x^i(t),y^a(t))$, reads \be (E_L):\qquad\frac{\xd x^i}{\xd t}=\sum_c\zr^i_c(x)y^c, \qquad
\frac{\xd}{\xd t}\left(\frac{\partial L}{\partial y^b}\right)= \sum_{a,c}{\mathcal C}^c_{ab}(x) y^a\frac{\partial
L}{\partial y^c}(x,y) + \sum_i\zs^i_b(x)\frac{\partial L}{\partial x^i}(x,y)\,.\label{EL2}
\ee
As one can see from (\ref{EL2}), the solutions are automatically
{\it admissible curves} in $E$, i.e. the velocity
$\bt(\zt\circ\zg)(t)$ is $\zr(\zg(t))$.

With this framework it is possible to write in a unified point of
view many equations of different mechanical systems that usually,
in the literature, appears as different ones (Classical
Euler-Lagrange equations, Lagrange-Poincar\'e equations after
reduction by the action of a Lie group, Euler-Poincar\'e
equations, etc.). (See \cite{CoLeMaMaMa,LMM} for applications of
the theory in the case when $(E,\ze)$ is  a Lie algebroid).

\section{Nonholonomic mechanics}\label{koi}
We can start from a general Lie algebroid $(E,\ze)$ keeping in
mind the standard case $E=\sT M$. A (linear) nonholonomic
Lagrangian system is determined by a Lagrangian function $L: E\ra
\R$ and a vector subbundle $D$, $\hbox{rank}\, D=n-r$, of the
bundle $E$. We will accept subbundles over a submanifold, so let
us denote $D_M=\zt(D)$. By $i_D: D\hookrightarrow E$ let us denote
the inclusion and by $i^*_D: E^*_{\mid D_M}\rightarrow D^*$ the
dual map.
%We have also the canonical map $\zp_D^*:\sT^*E_{\mid
%D}\ra \sT^* D=(\sT^*E_{\mid D})/(\sT D)^0$, where $(\sT
%D)^0\subset \sT^*E_{\mid D}$ is the annihilator of $\sT D$.

Because the solutions of the dynamics in $E$ should be admissible
curves, we need an {\it admissibility condition} ensuring that
there are admissible curves through every point of $D$. The
natural condition we take is $\zr(D)\subset \sT D_M$. By {\it
strong admissibility condition} we will mean that $D$ satisfies
the integrabilty condition with respect to both algebroid
structures: $\ze$ and $\ze^+$, i.e.
\be\label{ic}\zr(D)\subset \sT
D_M\quad\text{and}\quad\zs(D)\subset \sT D_M\,.
\ee
The Lagrangian function $L: E\rightarrow \R$ and the
vector subbundle $D$ define also the smooth map -- the {\it
constrained Tulczyjew differential} \be\label{nhc}\zL_L^D:
{D}\rightarrow \sT D^\ast,\qquad \zL_L^D=\sT
i_D^\ast\circ\ze\circ\xd L\ee covering the {\it constrained
Legendre map} \be\label{nh0}\zl_L^D:D\ra D^*,\qquad
\zl_L^D=i_{D}^\ast \circ \zl_L\,.
\ee
The diagram picture is the following
\be\label{nhc1} \xymatrix{ \sT^\ast E_{\mid D}\ar[rr]^{\ze}  && \sT (E^\ast_{\mid
D_M})\ar[rr]^{\sT{i_D^\ast}}&&
\sT D^\ast\ar[d]^{\zt_{D^\ast}}\\
D\ar[rr]^{\zl_L}\ar[u]^{\xd
L}\ar@{.>}[rrrru]^{\zL_L^D}\ar[rr]^{\zl_{L}}\ar@/d3ex/@{.>}[rrrr]_{\zl_L^D} &&
E^\ast_{\mid D_M}\ar[rr]^{i_D^\ast} && D^\ast}
\ee

\begin{definition}{\rm The bundle $D^\*$ we call
{\it the phase space} of the nonholonomic system and the implicit
differential equation ${\zG^D_L}=\zL_L^D(D)$, being a subset of
$\sT D^*$, we interpret as the {\it phase dynamics} of the
nonholonomic system. A curve $\gamma: I\rightarrow D$ is a
solution of the {\it nonholonomic Euler-Lagrange equation} if and
only if $\zL_L^D\circ \gamma$ is an admissible curve in $\sT
D^\*$, i.e. the nonholonomic Euler-Lagrange equation reads
\be\label{nEL}(nE_L):\qquad\qquad\bt(\zl^D_L\circ\zg)=\zL_L^D\circ \gamma\,.\ee}
\end{definition}

\bigskip To find the explicit form of the nonholonomic Euler-Lagrange equation (\ref{nEL}),
consider local coordinates $(x^I)=(x^i,x^\zi)$ on a open set ${U}$
of $M$ such that $D_M$ is
 determined by the constraint $x^\zi=0$.
A local basis $\{e_a\}_{a=1,\dots,n-r}$ of sections of $D$ we can
extend to local sections of $E$ and complete them to a local basis
of sections $\{e_a, e_{\alpha}\}$ of the vector bundle $E$. Then,
in coordinates $(x^I,y^A)=(x^i,x^\zi,y^a,y^\alpha)$ adapted to
this bases, the local equations defining the constrained subbundle
$D$ as a vector subbundle of $E$ over $D_M$ are $x^\zi=0$,
$y^\alpha=0$. Note that admissibility of the constraint $D$ means
that $\zr_a^\zi(x^i,0)=0$.

Taking local coordinates $(x^i, y^a)$ on $D$, we may write then
$i_D: D\hookrightarrow E$ as $i_D(x^i, y^a)=(x^i, 0, y^a, 0)$ and
$i_D^*(x^i,0,\zx_a,\zx_\za)=(x^i,\zx_a)$, so
$$\sT i_D^*(x^i,0,\zx_A,\dot{x}^j,0,\dot{\zx}_B)=(x^i,\zx_a,\dot{x}^j,\dot{\zx}_b)\,.$$
For the adapted local coordinates
$(x^i,x^\zi,y^a,y^\za,p_j,p_\zg,\zx_b,\zx_\zb)$ in $\sT^*E$, the
map $\ze$ reduced to $(\sT^*E)_{\mid D}$ takes values in
$\sT(E^*_{\mid D_M})$ (admissibility) and reads
$$\ze(x^i,0,y^a,0,p_J,\zx_B)=(x^i,0,\zx_A,\zr_d^j(x^i,0)y^d,0,{\mathcal
C}^C_{dB}(x^i,0)y^d\zx_C+\zs^J_B(x^i,0)p_J)\,.
$$
Therefore
\be\label{yy}\sT i_D^*\circ\ze(x^i,0,y^a,0,p_J,\zx_B)=(x^i,\zx_a,\zr_d^j(x^i,0)y^d,{\mathcal
C}^C_{db}(x^i,0)y^d\zx_C+\zs^J_b(x^i,0)p_J)
\ee
and
\begin{eqnarray*}
&&\zL_L^D(x^i, y^a)\\
&&=Ti_D^*(\ze(x^i, 0,y^a, 0, \frac{\partial L}{\partial
x^J}(x^i,0, y^a, 0), \frac{\partial L}{\partial
y^B}(x^i,0, y^a, 0))\\
&&=(x^i, \frac{\partial L}{\partial y^b}(x^i, 0,y^a, 0),
\rho^j_d(x^i,0) y^d, {{\mathcal C}}^C_{db}(x^i,0) y^d
\frac{\partial L}{\partial y^C}(x^i, 0,y^a,
0)+\sigma^J_b(x^i,0)\frac{\partial L}{\partial x^J}(x^i,0, y^a,
0))
\end{eqnarray*}
Therefore, locally,  the nonholonomic Euler-Lagrange equations read:
\begin{eqnarray}\label{nELe1}&\qquad
x^\zi=0,\quad y^\za=0,\quad \frac{d x^j}{dt}=\rho^j_d(x^i,0)y^d,\\
 &\frac{d}{dt}\left(\frac{\partial L}{\partial
y^b}(x^i, 0,y^a, 0)\right)= {{\mathcal C}}^C_{db}(x^i,0) y^d
\frac{\partial L}{\partial y^C}(x^i,0, y^a,
0)+\sigma^J_b(x^i,0)\frac{\partial L}{\partial x^J}(x^i, 0,y^a,
0))\,.\label{nELe2}
\end{eqnarray}
In the case of a Lie algebroid $\sigma^i_A= \rho^i_A$, and if the
subbundle $D$ is over the whole base manifold $M$, $D_M=M$, the
previous equations are precisely the nonholonomic equations
obtained in \cite{CoLeMaMa} (see Equations 3.8).

\section{The nonholonomic reduction of algebroids}

It has been recognized a long time ago \cite{Mr,SM} that
nonholonomic constraints may lead to certain {\it nonholonomic
brackets} that do not satisfy the Jacobi identity. We will show
that linear nonholonomic systems of mechanical type on general
algebroids are again systems on general algebroids. But even if we
start with a Lie algebroid, the new algebroid is, in general, no
longer a Lie algebroid but certain quasi-Lie algebroid associated
with a linear bi-vector field. This observation, made already in
\cite{GG} (see also \cite{LMMdD}), puts new light to the role of
{\it quasi-Poisson brackets}, i.e. the brackets represented by
arbitrary bi-vector fields and not satisfying, in general, the
Jacobi identity (see \cite{GG,LMMdD}). This shows, on the other
hand, that developing Mechanics on general algebroids makes sense,
as the reduction to a nonholonomic constraint will move us, in
general, from the Lie algebroid picture into a more general one.

Let $\ze$ be an algebroid structure on a vector bundle $E$ over
$M$ associated with the tensor $\zP_\ze$ with the local form
(\ref{tensor}). For a linear subbundle $D$ in $E$ over a
submanifold $D_M\subset M$, satisfying the strong admissibility
condition (\ref{ic}), consider a decomposition
\be\label{dec}E_{\mid D_M}=D\oplus_{D_M} D^\perp
\ee
and the associated projection $P:{E}_{\mid D_M}\ra{D}$. We can
construct the morphism ${\ze_P}: \sT^*D\rightarrow \sT D^*$ of
double vector bundles as follows.

To have the corresponding expressions in local coordinates,
consider local coordinates $(x^I)=(x^i,x^\zi)$ on a open set ${U}$
of $M$ such that $D_M$ is determined by the constraint $x^\zi=0$.
Local bases $\{e_a\}_{a=1,\dots,n-r}$ and
$\{e_\za\}_{\za=n-r+1,\dots,n}$ of sections of $D$ and $D^\perp$,
respectively, can be extended to a basis $e_A$ of local sections
of $E$. Then, we get the coordinates
$(x^I,y^A)=(x^i,x^\zi,y^a,y^\alpha)$ adapted to these bases and
the adapted coordinates
$(x^I,y^A,p_J,\zx_B)=(x^i,x^\zi,y^a,y^\alpha,p_j,p_\zg,\zx_b,\zx_\zb)$
in $\sT^\ast E$. We get automatically coordinates $(x^i,y^a)$ in
$D$ and $(x^i,y^a,y^\za)$ in $E_{\mid D_M}$, and the adapted
coordinates $(x^i,y^a,p_j,\zx_a)$ and $(x^i,y^A,p_j,\zx_B)$ in
$\sT^\ast D$ and $\sT^\ast(E_{\mid D_M})$, respectively.

Let us consider the phase lift $\sT^*P:\sT^*D\rel\sT^*(E_{\mid
D_M})$ which, as often happens which phase lifts, is not a map but
only a relation. In our local coordinates,
$$\sT^*P(x^i,y^a,p_j,\zx_b)=(x^i,y^A,p_j,\zx_b,0)\,.$$
We have also the embedding $i_{E_{\mid D_M}}:E_{\mid D_M}\ra E$
whose phase lift, restricted to $\sT^*E_{\mid D}$,
$$\sT^*(i_{E_{{\mid D}_{M}}}):\sT^*E_{\mid D}
\rel\sT^*(E_{\mid D_M})
$$
is the relation
$$\sT^*(i_{E_{{\mid
D}_{M}}})(x^i,0,y^a,0,p_J,\zx_B)=(x^i,y^A,p_j,\zx_B)\,.
$$
The composition of relations
$$\sT^*(i_{E_{{\mid D}_{M}}})^{-1}\circ\sT^*P:\sT^*D\rel\sT^*E_{\mid D}
$$
has the local form
$$\left(\sT^*(i_{E_{{\mid D}_{M}}})^{-1}\circ\sT^*P\right)(x^i,y^a,p_j,\zx_b)=
(x^i,0,y^a,0,p_J,\zx_b,0)\,.
$$
But, as $\zs_b^\zi(x^i,0)=0$ (strong admissibility condition), the
value of $\sT i_D^*\circ\ze$ in (\ref{yy}) does not depend on
$p_\zg$ and the composition
\be\label{nhp} {\ze_P}= \sT i^*_D\circ
\ze\circ \sT^*(i_{E_{{\mid D}_{M}}})^{-1}\circ\sT^*P: \sT^*D\ra
\sT D^* \ee
is a well-defined map which, in local coordinates,
reads
\be\label{yyy}{\ze_P}(x^i,y^a,p_j,\zx_b)=(x^i,\zx_a,\zr_d^j(x^i,0)y^d,{\mathcal
C}^c_{db}(x^i,0)y^d\zx_c+\zs^k_b(x^i,0)p_k)\,. \ee This is of
course an algebroid structure on the bundle $D$. This algebroid
structure can be described in a more straightforward way as
follows.

The decomposition (\ref{dec}) gives the dual decomposition
\be\label{dec1}E_{\mid D_M}^*=D^*\oplus_{D_M} (D^\perp)^* \ee and
the corresponding projection $P^*:E_{\mid D_M}^*\ra D^*$. This
projection does not give us a canonical projection of vectors and,
more generally, contravariant tensors on $E$ at points of $E_{\mid
D_M}^*$, unless they are tangent to $E_{\mid D_M}^*$. But the
tensor $\zP_\ze$ with the local form (\ref{tensor}) is a sum of
tensor product with at least one part in the product being
vertical. Moreover, due to strong admissibility condition for $D$,
there is a unique decomposition
$\zP_\ze=\zP_{\ze_D}+\zP_{\ze_{D^\perp}}$, where $\zP_{\ze_D}$ is
a linear tensor tangent to $D^*$ and the vertical parts of
$\zP_{\ze_{D^\perp}}$ are tangent to $(D^\perp)^*$, so they will
be killed by the projection, independently how the other part of
the product is. This gives a well-defined projection
$\sT(P^*)(\zP_\ze)=\zP_{\ze_P}$, when we restrict to the points of
$D^*$.

To put it differently, we can take any smooth projection
${\mathcal P}$ from a neighbourhood of $E_{\mid D_M}^*$ onto
$E_{\mid D_M}^*$. Then,
\be\label{gg}\sT(P^*)(\zP_\ze):=\sT(P^*\circ{\mathcal P})(\zP_\ze)
\ee
is a linear 2-contravariant tensor on $D$ which does not depend on
the choice of $\mathcal P$. This tensor defines on $D^*$ the {\it
nonholonomic bracket} $\{\cdot,\cdot\}_{\ze_P}$ associated with
the projection $P$.

On the level of the algebroid bracket $[\cdot,\cdot]_\ze$ this
procedure is the following. The strong admissibility condition for
$D$ implies that the bracket $[X,Y]_\ze$ of sections of $D$ is a
section of $E_{\mid D_M}$. Projecting this section to $D$ along
$P$ gives us a bracket
\be\label{nbra}[X,Y]_{\ze_P}=P[X,Y]_\ze
\ee
on sections of $D$ -- the {\it nonholonomic restriction of
$[\cdot,\cdot]$ along $P$}. This is an algebroid bracket with the
original anchors.

Of course, if the constraint bundle $D$ is a subalgebroid of $E$,
then it is a {\it holonomic} constraint and the projection $P$
plays no role: the tensor $\zP_\ze$ is tangent to $D$ and we just
take the restriction. It is however clear that, in the true
nonholonomic case, the nonholonomic bracket need not satisfy the
Jacobi identity, even when the original bracket does. In local
coordinates,
$$ \zP_{\ze_{P}} =\sum_{a,b,c}{\mathcal C}^c_{ab}(x^i,0)\zx_c \partial _{\zx_a}\otimes \partial _{\zx_b} +
\sum_{j,b}\zr^j_b(x^i,0) \partial _{\zx_b} \otimes \partial _{x^j}
- \sum_{a,j}\zs^j_a(x^i,0)\partial _{x^j} \otimes \partial
_{\zx_a}.
$$
and the corresponding nonholonomic bracket reads:
\begin{eqnarray*}
\{\zx_a, \zx_b\}_{\ze_{P}}&=&{\mathcal C}^c_{ab}(x^i,0)\zx_c, \\
\{\zx_b, x^j\}_{\ze_{P}}&=&\zr^j_b(x^i,0), \\
\{x^j, \zx_a\}_{\ze_{P}}&=&-\zs^j_a(x^i,0),\\
\{x^k, x^j\}_{\ze_{P}}&=& 0.
\end{eqnarray*}
In the case of a Lie algebroid, we have ${\mathcal
C}^c_{ab}=-{\mathcal C}^c_{ba}$, $\sigma^j_A=\rho^j_A$, and the
above bracket corresponds to the one introduced by \cite{SM} (see
also \cite{CoLeMaMa} and references therein) for the subbundles
over the total base $M$.

A natural question arises here: what Lie algebroid structures we
can obtain as nonholonomic restrictions of Lie algebroid brackets?
Of course, as the brackets must be skew-symmetric automatically,
the algebroids must be necessarily quasi-Lie algebroids. We will
call them {\it nonholonomic quasi-Lie algebroids}. In other words:
which linear bi-vector fields on a vector bundle $D \to M$ can be
obtained by projections of Poisson tensors from a bigger vector
bundle $E \to M$. For nonholonomic quasi-Lie algebras the answer
is simple.
\begin{theorem} Any quasi-Lie algebroid $D$ with the trivial anchor
is a nonholonomic quasi-Lie algebroid.
\end{theorem}
\begin{proof} As the anchor is trivial, we can clearly reduce to finite-dimensional
quasi-Lie algebras. Take a basis $e_i$ of $D$ and the
corresponding linear coordinates $y^i$, $i=1,\dots,n$. The
algebroid bracket in $D$ is then determined by the structure
constants $c_{ij}^k$,
$$[e_i,e_j]=c_{ij}^ke_k\,,$$
satisfying $c_{ij}^k=-c_{ji}^k$. Let us consider a new algebra $E$
with the basis $e_i,f_j$, $i=1,\dots,n$ and the bracket
$[\cdot,\cdot]_E$ for which $f_j$ are central elements and
$$[e_i,e_j]_E=c_{ij}^kf_k\,.$$
Since the algebra is 2-step nilpotent, it satisfies the Jacobi
identity and it is a Lie algebra. We can them embed $D$ in $E$ by
putting $\zf(e_i)=(e_i+f_i)$ and take the complementary subspace
$D^\perp$ as spanned by elements $e_i-f_i$. The projection $P:E\ra
D$ is therefore given by
$$P(e_i)=e_i+f_i,\quad P(f_i)=e_i+f_i\,.$$
The nonholonomic bracket on the embedded submanifold $D$ is
therefore
$$[e_i+f_i,e_j+f_j]=P(c_{ij}^kf_k)=c_{ij}^k(e_k+f_k)\,,$$
thus the original bracket in $D$.
\end{proof}

\medskip
For a general quasi-Lie algebroid the situation is much more
complicated and we do not know a full characterization of
nonholonomic quasi-Lie algebroids. Note however that not all
quasi-Lie algebroids are nonholonomic reductions of Lie
algebroids.
\begin{example}{\rm Take a quasi-Lie algebroid structure $\ze$ on
$D=\sT M$ with the anchor $\zr=\id_{\sT M}$. The corresponding
tensor on $\sT^* M$ has therefore the local form
$$\zP_\ze=\frac{1}{2}c_{ij}^k(x)p_k\partial_{p_i}\wedge\partial_{p_j}
+\partial_{p_i}\wedge\partial_{x^i}\,.
$$
Suppose that the corresponding bracket is a nonholonomic reduction
of a Lie algebroid bracket $[\cdot,\cdot]_E$ of a bigger vector
bundle $E$ along a projection $P$,
$$[X,Y]_\ze=P([X,Y]_E)\,.
$$
Since the anchor $\zr$ of $D$ is the restriction of the anchor
$\zr_E$ to $D$, the Lie algebroid $E$ satisfies $\zr_E(E_{\mid
M})=\sT M$, thus $E_{\mid M}=D\oplus_M K$, where $K$ is the kernel
of $\zr_E$ over $M$. For sections $X,Y$ of $D$ we have therefore
$\zr_E([X,Y]_E)=\zr([X,Y])$, thus, as the anchor map is for Lie
algebroids a homomorphism of brackets,
$$\zr([X,Y])=\zr_E([X,Y]_E)=[\zr_E(X),\zr_E(Y)]_{vf}
=[\zr(X),\zr(Y)]_{vf}\,,
$$
where $[\cdot,\cdot]_{vf}$ is the bracket of vector fields. We get
that the anchor $\zr$ maps the algebroid bracket into the bracket
of vector fields which is not the case for generally non-zero
structure functions $c_{ij}^k(x)$. }
\end{example}

\section{Nonholonomic systems of mechanical type}
\noindent Let us consider now on our general algebroid $(E,\ze)$ a
Lagrangian $L: E\rightarrow \R$ of mechanical type, i.e. \be\label
{ml} L(e)=\frac{1}{2}{G} (e, e)-V(\tau(e))\,, \ee where ${\mathcal
G}:E\times_M E\to \R$ is a bundle metric on $E$ and $V:
M\rightarrow \R$ is the potential function. Let us consider also a
a vector subbundle $D$ of $E$ over $D_M$ satisfying the strong
admissibility condition. Having the metric ${\mathcal G}$ in $E$,
we have the natural decomposition $E_{\mid
D_M}=D\oplus_{D_M}D^\perp$ with $D^\perp$ being the orthogonal
complement of $D$ with respect to the metric ${\mathcal G}$,
accompanied with the corresponding (this time -- orthogonal)
projection $P:E_{\mid D_M}\ra D$. The fundamental observation in
this case is the following.
\begin{theorem}\label{Theorem1}The nonholonomic Tulczyjew
differential ${\zL^D_L}$ coincides with the Tulczyjew differential
$\zL_l$ associated with the restricted Lagrangian $l=L_{\mid D}:
D\rightarrow \R$ and the nonholonomic algebroid structure $\ze_P$
on $D$. In other words, the nonholonomic Euler-Lagrange equations
associated with $L$ and the nonholonomic constraint $D$ in $E$
coincide with the non-constrained Euler-Lagrange equations for the
Lagrangian function $l$ on the algebroid $(D,{\ze_P})$.
\end{theorem}
\begin{proof}
Let us consider local coordinates $(x^I)=(x^i,x^\zi)$ on a open
set ${U}$ of $M$ such that $D_M$ is determined by the constraint
$x^\zi=0$. Take a local basis of orthonormal sections $\{e_a,
e_{\alpha}\}$ of $E$ adapted to the orthogonal decomposition
$E_{\mid D_M}=D\oplus_{\mid D_M} D^{\perp}$, i.e, $\hbox{span}
\{e_a\}=D$ and $\hbox{span} \{e_\alpha\}=D^\perp$. We get the
induced coordinates $(x^I,y^A)$ in $E$ in which the Lagrangian
takes the form
$$L(x,y)=\frac{1}{2}\sum _A(y^A)^2-V(x^i,x^\zi)\,.$$
Note that the `mass' is 1 in these coordinates. On $E_{\mid D_M}$
we have the induced coordinates $(x^i, y^a, y^{\alpha})$, so $D$
is locally determined by the constraints $x^\zi=0,y^{\alpha}=0$
and the restricted Lagrangian is:
\[
l(x^i, y^a)=\frac{1}{2}\sum_{a=1}^{r} (y^a)^2-V(x^i,0).
\]
In our coordinates, the projector ${P}$ reads
\[
{P}(x^i, y^a, y^{\alpha})=(x^i, y^a)
\]
and we deduce that $${\ze_P}(x^i,
y^a,p_j,\zx_b)=(x^i,\zx^a,\zr^j_d(x^i,0)y^d,{\mathcal
C}^c_{db}(x^i,0)y^d\zx_c+\zs^k_b(x^i,0)p_k)\,,
$$
so that
$$\ze_{P}\circ \xd l(x^i, y^a)
=\left(x^i, y^a, \rho^j_d(x^i,0)y^d, \sum_{c}{{\mathcal
C}}^c_{db}(x^i,0)y^d y^c-\sigma^k_b\frac{\partial V}{\partial
x^k}(x^i,0)\right)\,.
$$
We get therefore the Euler-Lagrange equations for the Lagrangian
$l$ on $D$ of the form
\begin{eqnarray*}
  \frac{d x^j}{dt}&=&\rho^j_a(x^i,0)y^a,\\
 \frac{d}{dt}\left(\frac{\partial l}{\partial
y^b}(x^i,y^a)\right)&=& {{\mathcal C}}^c_{db}(x^i,0) y^d
y^c+\sigma^j_b(x)\frac{\partial l}{\partial x^j}(x^i,y^a),
\end{eqnarray*}
which are exactly the nonholonomic equations (\ref{nELe1}) and
(\ref{nELe2}), if we take to account that
$$\frac{\partial L}{\partial y^C}(x^i,0, y^a, 0)=\zd^C_ay^a$$
and
$$\sigma^\zg_b(x^i,0)=0\,.$$
\end{proof}
The above theorem makes it clear that, at least for Lagrangians of
mechanical type, the nonholonomic Euler-Lagrange equations are
just Euler-Lagrange equations, but on the nonholonomic restriction
of the algebroid. In any case, passing to a liner nonholonomic
constraint does not move us out of the Mechanics on algebroids.
Our theory is therefore complete with respect to passing to the
nonholonomic case, what was the main problem in understanding the
nonholonomic constraint in the Lie algebroid (thus canonical)
case.

\medskip
Having interpreted the nonholonomic Euler-Lagrange equations as
just Euler-Lagrange equations, but on a reduced algebroid, we can
use the generalization of Noether Theorem formulated for general
algebroids in \cite{GGU} to obtain a nonholonomic Noether Theorem,
at least for Lagrangians of mechanical type.

The generalized Noether Theorem for a general algebroid $(E,\ze)$
is based on the concept of the complete lift $\dt^\ze(K)$ of
tensor fields $K$ being sections of the tensor products
$E^{\otimes k}$ to the corresponding contravariant tensor fields
on $E$. This is a natural generalization of the standard concept
of the tangent lift $\dt$ which lifts contravariant tensors on a
manifold $M$ to the corresponding tensor fields on $\sT M$, (cf.
\cite{GU1,YI}).

For a vector bundle $\zt:E\ra M$, let $\otimes ^k(\zt)$ be the
space of sections of the tensor-product bundle $E^{\otimes k}$
over $M$. With any tensor field $K\in \otimes ^k(\zt)$ we can
associate the linear function $\zi(K)$ on the dual bundle
$(E^{\otimes k})^*=(E^*)^{\otimes k}$ and the vertical lift
$\textsf{v}_\zt(K) \in \otimes ^k(\ezt_E)$  In local coordinates,
$$ \zi(f^{a_1\dotsi a_k}(x) e_{a_1}\otimes \cdots \otimes
e_{a_k}) = f^{a_1\dotsi a_k}(x)\zx_{a_1\dotsi a_k}
$$
and
$$ \textsf{v}_\zt(f^{a_1\dotsi a_k}(x) e_{a_1}\otimes \cdots \otimes
e_{a_k}) = f^{a_1\dotsi a_k}(x) \partial _{y^{a_1}} \otimes \cdots
\otimes \partial _{y^{a_k}}\,.
$$
A particular case of the vertical
lift is the lift $\vt(K)$ of a contravariant tensor field $K$ on
$M$ into a contravariant tensor field on $\sT M$. It is well known
(see \cite{YI,GU1}) that in the case of $E=\sT M$ we have also the
tangent lift $\dt \colon \otimes (\zt_M) \rightarrow \otimes
(\zt_{\ss T M})$ which is a $\vt$-derivation. What has been done
for the tangent bundle, can be repeated in the case of an
arbitrary algebroid $(E,\ze)$. Note first that we can extend $\ze$
naturally to mappings (cf. \cite{GU3, GU2})
$$ \ze^{\otimes r} \colon \otimes ^r_E \sT^\* E \longrightarrow
\sT \otimes ^r_M E^\*,\quad r\ge 0.
$$

\begin{theorem}{\em \cite{GU3, GU2}}
Let $(E,\ze)$ be an algebroid. For $K\in \otimes ^k(\zt)$, $k\ge
0$, the equality
\be \zi(\dt^\ze(K)) = \dt (\zi(K))\circ \ze^{\otimes k} \ee
defines the  tensor field $\dt^\ze(K) \in \otimes ^k(\zt_E)$ which
is linear and the mapping
$$ \dt^\ze \colon \oplus_k\otimes^k (\zt) \longrightarrow \oplus_k\otimes^k (\zt_E) $$
is a $\vt$-derivation of degree $0$. In local coordinates, the
lifts of functions on $M$ and sections of $E$ read:
\bea
\dt^\ze(f(x))& =&  y^a \zr^i_a(x) \frac{\partial f}{\partial
x^i}(x),\\
\dt^\ze(f^a(x)e_a) &=& f^a(x)\zs^i_a(x)\partial _{x^i} + \left(
y^a \zr^i_a(x) \frac{\partial f^c}{\partial x^i}(x) + {\mathcal
C}^c_{ab}(x)y^af^b(x) \right) \partial _{y^c}\,.
\eea
Conversely, if $D\colon \oplus_k\otimes^k (\zt) \longrightarrow
\oplus_k\otimes^k (\zt_E)$ is a $\vt$-derivation of degree $0$
such that $D(K)$ is linear for each $K \in \otimes ^1(\zt)$, then
there is an algebroid structure $\ze$ on $\zt\colon E\rightarrow
M$ such that $D = \dt^\ze$. This algebroid structure is a Lie
algebroid if and only if
$$\dt^\ze([X,Y]_\ze) = [\dt^\ze(X), \dt^\ze(Y)]_{vf}
$$
for all $X,Y\in\otimes ^1(\zt)$.
\end{theorem}
Suppose now that we are dealing with a Lagrangian $L$ of
mechanical type and a nonholonomic constraint subbundle $D$ of
$(E,\ze)$ satisfying the strong admissibility assumption, so that
we have also the nonoholomic reduced algebroid $(D,\ze_P)$. We
will say that a pair $(X,f)$ consisting of a section $X$ of $D$
and a function $f$ on $D_M$ is a {\it symmetry of the nonholonomic
problem} associated with the Lagrangian $L$, if
\be\label{Ne}
{\dt^{\ze_P}(X)}(l)=\dt^{\ze_{P}}(f)\,,
\ee
where $l=L_{\mid D}$ is the restriction of the Lagrangian to the
constraint.
\begin{theorem}{[Nonholonomic Noether Theorem]}
The following are equivalent:
\begin{itemize}
\item[(a)] the pair $(X,f)$ is a symmetry of the nonholonomic
problem associated with the Lagrangian $l$;

\item[(b)] the function $\zi_{X}-f\circ\zt_{D^*}$ on $D^*$ is a
constant of the motion of the nonholonomic phase dynamics on
$D^*$, i.e. $\dt(\zi_{X}-f\circ\zt_{D^*})$ vanishes on
${\zG^D_L}$;

\item[(c)] the function $$(\zi_{X}-f\circ\zt_{D^*})\circ\zl_L^D
$$
is a constant of the motion for the nonholonomic Euler-Lagrange
equation.
\end{itemize}
\end{theorem}
\begin{proof}
The proof is completely analogous to that for  Theorem 4 in
\cite{GGU}.
\end{proof}
For a section $X$ of $D$,  not verifying necessarily condition (\ref{Ne}), we obtain
\be\label{nme}
{\dt^{\ze_P}(X)}(l)\circ\zg(t)=\frac{\xd}{\xd t}\left(
\zi_{X}\circ\zl_L^D\circ\zg(t)\right)
\ee
for any solution $\zg:\R\ra D$ of the nonholonomic Euler-Lagrange
equation. This last equation can be interpreted as a general
version of the {\it nonholonomic momentum equation} studied for
several authors (see \cite{CoLeMaMa,BKMM} and references therein)
for nonholonomic systems with symmetry.

\section{Examples}

\begin{example}{The  Chaplygin sleigh.}

{\rm As an example of nonholonomic system on a Lie algebra, we
study   the  Chaplygin sleigh which describes a rigid body sliding
on a plane. The body is supported in three points, two of which
slides freely without friction while the third point is a knife
edge.  This imposes the constraint of no motion orthogonal to this
edge (see \cite{Cha,NeFu}).

The configuration space before reduction is the Lie group
$G=SE(2)$ of the Euclidean motions of the 2-dimensional plane
$\R^2$. We will need in the sequel to fix some notation about the
Lie algebra $\mathfrak{se}(2)$. First of all its elements are
matrices of the form
\[
\hat{\xi}=\left(
\begin{array}{ccc}
0&\xi_3&\xi_1\\
-\xi_3&0&\xi_2\\
0&0&0
\end{array}
\right)
\]
and a basis of the Lie algebra $\mathfrak{se}(2)\cong \R ^3$
is given by
\[
E_1=\left(
\begin{array}{ccc}
0&0&1\\
0&0&0\\
0&0&0
\end{array}
\right) ,\qquad E_2=\left(
\begin{array}{ccc}
0&0&0\\
0&0&1\\
0&0&0
\end{array}
\right) ,\quad
E_3=\left(
\begin{array}{ccc}
0&-1&0\\
1&0&0\\
0&0&0
\end{array}
\right).
\]
We have that
\[
{}[E_3,E_1]=E_2,\quad [E_2, E_3]=E_1,\quad [E_1, E_2]=0.
\]
An element $\xi \in \mathfrak{se}(2)$ is of the form
\[
\xi =v_1\, E_1+v_2\, E_2+\omega \, E_3.
\]

The Chaplygin system is described by the kinetic Lagrangian function
\[
\begin{array}{rrcl}
L:&\mathfrak{se}(2)&\longrightarrow&\R\\
  &(v_1, v_2, \omega)&\longmapsto&\frac{1}{2}\left[ (J+m(a^2+b^2))\omega^2 + mv_1^2+m v_2^2-2bm\omega v_1-2am\omega v_2\right]
  \end{array}
\]
where   $m$ and $J$ denotes the mass and moment of inertia of the
sleigh relative to the contact point and $(a, b)$ represents the
position of the center of mass with respect to the body frame
determined placing the origin at the contact point and the first
coordinate axis in the direction of the knife axis. Additionally,
the   system is subjected to the nonholonomic constraint
determined by the linear subspace of $\mathfrak{se}(2)$:
\[
D=\{(v_1, v_2, \omega)\in se(2)\; |\; v_2=0\}\, .
\]
Instead of $\{E_1, E_2, E_3\}$ we take the basis of $\mathfrak{se}(2)$:
\[
\{e_1=E_3, e_2=E_1, e_3= -ma E_3-mab E_1+(J+ma^2) E_2\}
\]
which is  a basis adapted to the decomposition $D\oplus D^\perp$; $D=\hbox{span }\{ e_1, e_2\}$ and $D^\perp=\hbox{span }\{ e_3\}$.

In the induced coordinates $(y^1, y^2)$ on $D$ the restricted lagrangian is
\[
l(y^1, y^2)=\frac{1}{2}\left[ (J+m(a^2+b^2))(y^1)^2 +
m(y^2)^2-2bmy^1y^2\right]\; ,
\]
and moreover,
\[
{}[e_1, e_2]_{\ze_P}=\frac{ma}{J+ma^2} e_1+\frac{mab}{J+ma^2}e_2\;
,\] Therefore, ${\mathcal C}^1_{12}=\frac{ma}{J+ma^2}$ and
${\mathcal C}^2_{12}=\frac{mab}{J+ma^2}$.

Then,
 $$\ze_P\circ \xd l(y^1, y^2))
=\left( (J+m(a^2+b^2))y^1 -bm y^2, my^2-bmy^1, -ma y^1y^2,
ma(y^1)^2\right).
$$
In consequence, the equations of motion are:
\begin{eqnarray*}
(J+m(a^2+b^2))\dot{y}^1 -bm \dot{y}^2&=& -ma y^1y^2\\
m\dot{y}^2-bm\dot{y}^1&=&ma(y^1)^2
\end{eqnarray*}
or,
\begin{eqnarray*}
\dot{y}^1 &=& \frac{ma}{J+ma^2}y^1\left(by^1-y^2\right)\\
\dot{y}^2&=&
\frac{ma}{J+ma^2}y^1\left((J+m(a^2+b^2))y^1-by^2\right).
\end{eqnarray*}

}

\end{example}

\begin{example}{The snakeboard}
{\rm

 As a mechanical system the snakeboard has as configuration
space $Q=SE(2)\times T^2$ with coordinates $(x, y, \theta,
\psi, \phi)$ (see ~\cite{BuZe,LMMdD,MaKo}).

\medskip
\begin{center}
\includegraphics{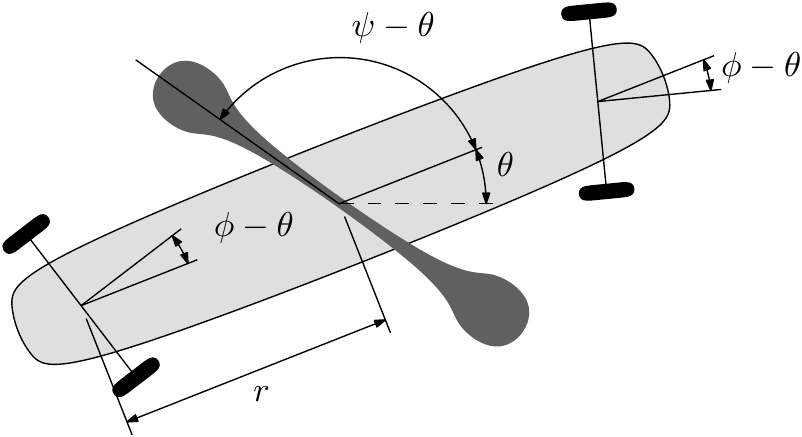}
\end{center}

\medskip
The  nonholonomic dynamics is described by
\begin{itemize}
\item The \emph{Lagrangian}
\begin{eqnarray*}
 L(q, \dot{q})= \frac{1}{2}m (\dot{x}^2+\dot{y}^2) + \frac{1}{2}
 (J+2J_1)\dot{\theta}^2
 +\frac{1}{2}J_0(\dot{\theta}+\dot{\psi})^2
 +J_1\dot{\phi}^2,
 \end{eqnarray*}
where $m$ is the total mass of the board, $J>0$ is the moment of
inertia of the board, $J_0>0$ is the moment of inertia of the
rotor of the snakeboard mounted on the body's center of mass and
$J_1>0$  is the moment of inertia of each wheel axles. The
distance between the center of the board and the wheels is denoted
by $r$. For simplicity (see  \cite{MaKo}), we assume that $J+J_0+2J_1=mr^2$.

\item The \emph{nonholonomic constraints} induced by the non sliding condition in the
sideways direction of the wheels:
\begin{align*}
-\dot{x}\sin (\theta+\phi) +\dot{y}\cos (\theta+\phi)
-r\dot{\theta}\cos\phi&=0\\
-\dot{x}\sin (\theta-\phi) +\dot{y}\cos (\theta-\phi)
+r\dot{\theta}\cos\phi&=0.
\end{align*}
\end{itemize}

Observe that the Lagrangian is induced by the riemannian  metric ${\mathcal G}$ on $Q$,
\[
{\mathcal G}=mdx^2+mdy^2+mr^2d\theta^2+J_0 d\theta\otimes\psi+J_0 d\psi\otimes
d\theta + J_0d\psi^2+2J_1 d\phi^2.
\]

 The constraint  subbundle $\tau_D: D\longmapsto Q$ is
 \[
{D}=\text{span}\left\{ e_1=\frac{\partial}{\partial
\psi}, e_2=\frac{\partial}{\partial \phi},
e_3=a\frac{\partial}{\partial x}+b\frac{\partial}{\partial
y}+c\frac{\partial}{\partial \theta}\right\}.
\]
where
\begin{eqnarray*}
a&=&-r(\cos\phi\cos(\theta-\phi)+\cos\phi\cos(\theta+\phi))=-2r\cos^2\phi\cos\theta\\
b&=&-r(\cos\phi\sin(\theta-\phi)+\cos\phi\sin(\theta+\phi))=-2r\cos^2\phi\sin\theta\\
c&=&\sin(2\phi).
\end{eqnarray*}

The orthogonal complement of $D$ is spanned by
\[
{D}^\perp=\text{span}\left\{ e_4=-b\frac{\partial}{\partial x}+a\frac{\partial}{\partial y},
e_5=(cJ-cmr^2)\frac{\partial}{\partial x}
+am\frac{\partial}{\partial \theta}-am\frac{\partial}{\partial \psi}\right\}.
\]
In the induced coordinates $(x, y, \theta,
\psi, \phi, y^1, y^2, y^3)$ on $D$  the restricted lagrangian is
\[
l((x, y, \theta,
\psi, \phi, y^1, y^2, y^3)=2mr^2\cos^2\phi (y^3)^2+J_0 c y^1y^3+\frac{1}{2}J_0 (y^1)^2+J_1 (y^2)^2.
\]
where now the nonholonomic constraints are rewritten as: $y^4=0$ and $y^5=0$.
After some straightforward computations we deduce that
\begin{eqnarray*}
{}[e_1, e_2]_{\ze_P}&=&0,
\\
{}[e_1, e_3]_{\ze_P}&=&0,\\
{}[e_2, e_3]_{\ze_P}&=&
\frac{2mr^2\cos^2\phi}{mr^2-J_0\sin^2\phi} e_1
-\frac{(mr^2+\cos 2\phi)\tan\phi}{mr^2-J_0\sin^2\phi} e_3.
\end{eqnarray*}
Thus, the unique non vanishing structure functions are:
\[
{\mathcal C}_{23}^1=-{\mathcal C}_{32}^1=\frac{2mr^2\cos^2\phi}{mr^2-J_0\sin^2\phi}, \qquad
{\mathcal C}_{23}^3=-{\mathcal C}_{32}^3=-\frac{(mr^2+\cos 2\phi)\tan\phi}{mr^2-J_0\sin^2\phi}  \; .
\]

%Finally,  $$\ze_P\circ \xd l
%  (x, y, \theta,
%\psi, \phi; y^1, y^2, y^3)
%=
%\left(
%\begin{array}{c}
%x, y, \theta,
%\psi, \phi {\bf ;}\;
%\zx_1=J_0(y^1+cy^3), \zx_2=2J_1y^2, \zx_3=4mr^2\cos^2\phi y^3+J_0 cy^1{\bf ;}\\
%ay^3, by^3, cy^3, y^1, y^2 {\bf ;}\; 0, 0,
%2J_0\cos(2\phi)y^1y^2-2mr^2\sin(2\phi)y^2y^3
%\end{array}
%\right).
%$$
Therefore, the equation of motion of the snakeboard are
\begin{eqnarray*}
&\dot{x}=ay^3,\qquad \dot{y}=by^3, \qquad \dot{\theta}=cy^3,\qquad \dot{\psi}=y^1, \qquad \dot{\phi}=y^2,&\\
&\frac{d}{dt}(y^1+cy^3)=0,\qquad \frac{d}{dt}(y^2)=0,&\\
&\frac{d}{dt}(4mr^2\cos^2\phi y^3+J_0 cy^1)= 2J_0\cos(2\phi)y^1y^2-2mr^2\sin(2\phi)y^2y^3&\; .
\end{eqnarray*}

}
\end{example}

\end{document}